\input amstex
\input amsppt.sty   
\hsize 14.5cm
\vsize 19.0cm
\def\nmb#1#2{#2}         
\def\totoc{}             
\def\idx{}               
\def\ign#1{}             

\redefine\o{\circ}

\define\al{\alpha}
\define\be{\beta}
\define\ga{\gamma}
\define\de{\delta}

\define\la{\lambda}

\define\ta{\tau}
\define\ph{\varphi}

\predefine\ii{\i}
\redefine\i{^{-1}}
\define\row#1#2#3{#1_{#2},\ldots,#1_{#3}}
\define\x{\times}
\define\Imm{\operatorname{Imm}}
\define\Diff{\operatorname{Diff}}
\define\Emb{\operatorname{Emb}}
\define\oin{\operatorname{int}}
\def\today{\ifcase\month\or
 January\or February\or March\or April\or May\or June\or
 July\or August\or September\or October\or November\or December\fi
 \space\number\day, \number\year}
\topmatter
\title The Action of the Diffeomorphism Group \\
on the Space of Immersions \endtitle
\author  Vicente Cervera\\
Francisca Mascar\'o\\
Peter W. Michor  \endauthor
\affil
Departamento de Geometr\'\ii a y Topolog\'\ii a\\
Facultad de Matem\'aticas\\
Universidad de Valencia\\ \\
Institut f\"ur Mathematik, Universit\"at Wien,\\
Strudlhofgasse 4, A-1090 Wien, Austria.
\endaffil
\address{Departamento de Geometr\'\ii a y Topolog\'\ii a,
Facultad de Matem\'aticas,
Universidad de Valencia,
E-46100 BURJASSOT,
VALENCIA,  SPAIN}\endaddress
\email mascaro\@evalun11.bitnet \endemail
\address{Institut f\"ur Mathematik, Universit\"at Wien,
Strudlhofgasse 4, A-1090 Wien, Austria}\endaddress
\date{1989}\enddate
\thanks{This paper was prepared during a stay of the third
author in Valencia, by a grant given by Coseller\'\ii a de Cultura,
Educaci\'on y Ciencia, Generalidad Valenciana.
The two first two authors were partially supported by the CICYT grant
n. PS87-0115-G03-01.}\endthanks
\abstract{ We study the action of the 
diffeomorphism group $\Diff(M)$ on the space of proper immersions 
$\Imm_{\text{prop}}(M,N)$ by composition from the right. 
We show that smooth transversal slices exist through 
each orbit, that the quotient space is Hausdorff and is stratified 
into smooth manifolds, one for each conjugacy class of isotropy 
groups.}\endabstract
\endtopmatter

\leftheadtext{\smc Cervera, Mascar\'o, Michor}
\rightheadtext{\smc The space of immersions}
\document

\heading Table of contents \endheading
\noindent Introduction \par
\noindent 1. Regular orbits\par
\noindent 2. Some orbit spaces are Hausdorff\par
\noindent 3. Singular orbits\par

\heading\totoc Introduction \endheading

Let $M$ and $N$ be smooth finite dimensional manifolds, connected and 
second countable without boundary such that $\dim M\leq\dim N$.
Let $\Imm(M,N)$ be the set of all immersions from $M$ into $N$.
It is an open subset of the smooth manifold $C^\infty(M,N)$, see
our main reference \cite{Michor, 1980c}, so it is itself a smooth
manifold. We also consider the smooth Lie group $\Diff(M)$ of
all diffeomorphisms of $M$. We have the canonical right action
of $\Diff(M)$ on $\Imm(M,N)$ by composition.

The space $\Emb(M,N)$ of embeddings from $M$ into $N$ is an
open submanifold of $\Imm(M,N)$ which is stable under the right
action of the diffeomorphism group. Then $\Emb(M,N)$ is the
total space of a smooth principal fiber bundle with structure
group the diffeomorphism group; the base is called $B(M,N)$, it
is a Hausdorff smooth manifold modeled on nuclear (LF)-spaces.
It can be thought of as the "nonlinear Grassmannian" of all
submanifolds of $N$ which are of type $M$. This result is based
on an idea implicitly contained in \cite{Weinstein, 1971}, it
was fully proved by \cite{Binz-Fischer, 1981} for compact $M$
and for general $M$ by \cite{Michor, 1980b}. The clearest
presentation is in \cite{Michor, 1980c, section 13}.
If we take a Hilbert space $H$ instead of $N$, then  $B(M,H)$ is
the classifying space for $\Diff(M)$ if $M$ is compact, and the
classifying bundle $\Emb(M,H)$ carries also a universal
connection. This is shown in \cite{Michor, 1988}.

The purpose of this note is to present a generalization of this
result to the space of immersions. It fails in general, since
the action of the diffeomorphism group is not free. Also we were not
able to show that the orbit space $\Imm(M,N)/\Diff(M)$ is Hausdorff. 
Let $\Imm_{\text{prop}}(M,N)$ be the space of all proper immersions.
Then $\Imm_{\text{prop}}(M,N)/\Diff(M)$
turns out to be Hausdorff, and the
space of those immersions, on which the diffeomorphism group
acts free, is open and is the total space of a smooth principal
bundle with structure group $\Diff(M)$ and a smooth manifold as
base space. For the immersions on which $\Diff(M)$ does not act
free we give a slice theorem which is explicit enough to
describe the stratification of the orbit space in detail.
The results are new and interesting even in the special case of
the loop space $C^\infty(S^1,N)\supset \Imm(S^1,N)$.

The main reference for manifolds of mappings is \cite{Michor,
1980c}. But the differential calculus used there is a little
old fashioned now, so it should be supplemented by the convenient
setting for differential calculus presented in
\cite{Fr\"olicher-Kriegl, 1988}.

If we assume that $M$ and $N$ are real analytic manifolds with $M$ 
compact, then all infinite dimensional spaces become real analytic 
manifolds and all results of this paper remain true, by applying 
the setting of \cite{Kriegl-Michor, 1990}.

\heading\totoc\nmb0{1}. Regular orbits \endheading

\subheading{\nmb.{1.1}. Setup} 
Let $M$ and $N$ be smooth finite dimensional manifolds,
connected and second
countable without boundary, and suppose that $\dim M\leq\dim N$.
Let $\Imm(M,N)$ be the manifold of all immersions from $M$ into $N$
and let $\Imm_{\text{prop}}(M,N)$ be the open submanifold of all proper
immersions.

Fix an immersion $i$. We will now describe some data for $i$
which we will use throughout the paper. If we need these data for
several immersions, we will distinguish them by appropriate 
superscripts.

First there are sets
$W_\al\subset\overline{W}_\al\subset U_\al\subset M$ such that
$(W_\al)$ is an open cover of $M$, $\overline W_\al$ is compact,
and $U_\al$ is an open locally finite cover of $M$, each $W_\al$
and $U_\al$ is connected, and such that $i|U_\al:U_\al\to N$ is
an embedding for each $\al$.

Let $g$ be a fixed Riemannian metric on $N$ and let $\exp^N$ be
its exponential mapping. Then let 
$p:\Cal N(i)\to M$ be the {\it normal bundle} of $i$, defined in
the following way: For $x \in M$ let $\Cal N(i)_x:=
(T_xi(T_xM))^\bot\subset T_{i(x)}N$ be the $g$-orthogonal
complement in $T_{i(x)}N$. Then 
$$\CD
\Cal N(i) @>\bar i>> TN\\
@VpVV                @VV\pi_NV\\
M         @>>i>      N
\endCD$$
is a vector bundle homomorphism over $i$, which is fiberwise
injective.  

Now let $U^i=U$ be an open neighborhood of the zero section which
is so small that 
$(\exp^N\o\bar i)|(U|U_\al): U|U_\al \to N$
is a diffeomorphism onto its image which describes a tubular
neighborhood  of the submanifold $i(U_\al)$ for each $\al$.
Let 
$$\ta=\tau^i:= (\exp^N\o\bar i)|U:\Cal N(i)\supset U\to N.$$
It will serve us as a substitute for a tubular neighborhood of $i(M)$.

\subheading{\nmb.{1.2}. Definition} An immersion $i\in
\Imm(M,N)$ is called {\it free }if $\Diff(M)$ acts freely on it,
i.e. if $i\o f=i$ for $f\in \Diff(M)$ implies $f=Id_M$.
Let $\Imm_{\text{free}}(M,N)$ denote the set of all free immersions.

\proclaim{\nmb.{1.3}. Lemma} Let $i\in \Imm(M,N)$ and let 
$f\in\Diff(M)$ have a fixed point $x_0\in M$ and satisfy 
$i\o f=i$. Then $f=Id_M$.
\endproclaim
\demo{Proof} We consider the sets $(U_\al)$ for the immersion $i$ 
of \nmb!{1.1}. Let us investigate $f(U_\al)\cap U_\al$.
If there is an $x\in U_\al$ with $y=f(x)\in U_\al$, we have
$(i|U_\al)(x)=((i\o f)|U_\al)(x)=(i|U_\al)(f(x))=(i|U_\al)(y)$.
Since $i|U_\al$ is injective we have $x=y$, and  
$$f(U_\al)\cap U_\al=\{x\in U_\al: f(x)=x\}.$$
Thus $f(U_\al)\cap U_\al$ is closed in $U_\al$. Since it is also
open and since $U_\al$ is connected, we have 
$f(U_\al)\cap U_\al=\emptyset$ or $=U_\al$.

Now we consider the set $\{x\in M: f(x)=x\}$. We have just shown
that it is open in $M$. Since it is also closed and contains the
fixed point $x_0$, it coincides with $M$.
\qed\enddemo

\proclaim{\nmb.{1.4}. Lemma} If for an immersion $i\in
\Imm(M,N)$ there is a point in $i(M)$ with only one preimage,
then $i$ is a free immersion.
\endproclaim
\demo{Proof} Let $x_0\in M$ be such that $i(x_0)$ has only one
preimage. If $i\o f=i$ for $f\in \Diff(M)$ then $f(x_0)=x_0$ and
$f=Id_M$ by lemma \nmb!{1.3}.
\qed\enddemo

Note that there are free immersions without a point in $i(M)$ with 
only one preimage: Consider a figure eight which consists of two 
touching circles. Now we may map the circle to the figure eight by 
going first three times around the upper circle, then twice around the 
lower one. This immersion $S^1\to \Bbb R^2$ is free.

\proclaim{\nmb.{1.5}. Theorem} Let $i$ be a free immersion 
$M\to N$. Then there is an open neighborhood $\Cal W(i)$ in $\Imm(M,N)$ 
which is saturated for the $\Diff(M)$-action and which splits 
smoothly as 
$$\Cal W(i)=\Cal Q(i)\x \Diff(M).$$
Here $\Cal Q(i)$ is a smooth splitting submanifold of $\Imm(M,N)$, 
diffeomorphic to an open neighborhood of 0 in $C^\infty(\Cal N(i))$. 
In particular the space $\Imm_{\text{free}}(M,N)$ 
is open in $C^\infty(M,N)$.

Let $\pi:\Imm(M,N)\to \Imm(M,N)/\Diff(M)=B(M,N)$ be the projection onto the 
orbit space, which we equip with the quotient topology. Then 
$\pi| \Cal Q(i):\Cal Q(i)\to \pi(\Cal Q(i))$ is bijective 
onto an open subset of the quotient. If $i$ runs through 
$\Imm_{\text{free,prop}}(M,N)$ of all free and proper immersions 
these mappings define a smooth atlas for 
the quotient space, so that 
$$(\Imm_{\text{free,prop}}(M,N), \pi, \Imm_{\text{free,prop}}(M,N)/\Diff(M), 
\Diff(M))$$
is a smooth principal fiber bundle with structure group $\Diff(M)$.
\endproclaim

The restriction to proper immersions is necessary because we are only
able to show that $\Imm_{\text{prop}}(M,N)/\Diff(M)$ is Hausdorff in
section \nmb!{2} below.

\demo{Proof}
We consider the setup \nmb!{1.1} for the free immersion $i$. Let 
$$\Cal U(i):=\{j\in \Imm(M,N):j(\overline{W}^i_\al)\subseteq 
     \tau^i(U^i| U^i_\al)\text{ for all }\al, j \sim i \},$$
where $j\sim i$ means that $j=i$ off some compact set in $M$. 
Then by \cite{Michor, 1980c, section 4} the set $\Cal U(i)$ is an 
open neighborhood of $i$ in $\Imm(M,N)$.
For each $j\in\Cal U(i)$ we define 
$$\align 
&\ph_i(j):M\to U^i\subseteq \Cal N(i),\\
&\ph_i(j)(x):=(\tau^i| (U^i| U^i_\al))\i(j(x))
     \text{ if } x\in W^i_\al.
\endalign$$
Then $\ph_i:\Cal U(i)\to C^\infty(M,\Cal N(i))$ is a mapping which is 
bijective onto the open set 
$$\Cal V(i):=\{h\in C^\infty(M,\Cal N(i)): 
     h(\overline W^i_\al)\subseteq U^i| U^i_\al
     \text{ for all }\al, h\sim 0\}$$
in $C^\infty(M,\Cal N(i))$. Its inverse is given by the smooth 
mapping $\tau^i_*:h\mapsto \tau^i\o h$, see \cite{Michor, 1980c, 
10.14}. We claim that $\ph_i$ is itself a smooth mapping: recall the 
fixed Riemannian metric $g$ on $N$; $\tau^i$ is a local 
diffeomorphism $U^i\to N$, so we choose the exponential mapping with 
respect to $(\tau^i)^*g$ on $U^i$ and that with respect to $g$ on 
$N$; then in the canonical chart of $C^\infty(M,U^i)$ centered at 0 
and of $C^\infty(M,N)$ centered at $i$ as described in \cite{Michor, 
1980c, 10.4}, the mapping $\ph_i$ is just the identity.

We have $\tau^i_*(h\o f)=\tau^i_*(h)\o f$ for those 
$f\in\Diff(M)$ which are near enough to the identity so that 
$h\o f\in \Cal V(i)$. We consider now the open set 
$$\{h\o f: h\in \Cal V(i), 
     f\in \Diff(M)\}\subseteq C^\infty((M,U^i)).$$
Obviously we have a smooth mapping from it into 
$C^\infty_c(U^i)\x \Diff(M)$ given by 
$h\mapsto (h\o(p\o h)\i,p\o h)$, where $C^\infty_c(U^i)$ is the space 
of sections with compact support of $U^i\to M$.
So if we let 
$\Cal Q(i):= \tau^i_*(C^\infty_c(U^i)\cap\Cal V(i))\subset \Imm(M,N)$
we have 
$$\Cal W(i):= \Cal U(i)\o \Diff(M)\cong \Cal Q(i)\x \Diff(M) \cong
     (C^\infty_c(U^i)\cap\Cal V(i))\x \Diff(M),$$
since the action of $\Diff(M)$ on $i$ is free.
Consequently $\Diff(M)$ acts freely on each immersion in $\Cal W(i)$, 
so $\Imm_{\text{free}}(M,N)$ is open in $C^\infty(M,N)$. Furthermore 
$$\pi|\Cal Q(i): 
\Cal Q(i)\to \Imm_{\text{free}}(M,N)/\Diff(M)$$
is bijective onto an open set in the quotient.

We now consider 
$\ph_i\o(\pi|\Cal Q(i))\i:\pi(\Cal Q(i))\to C^\infty(U^i)$
as a chart for the quotient space. In order to investigate the chart 
change let $j\in \Imm_{\text{free}}(M,N)$ be such that 
$\pi(\Cal Q(i))\cap\pi(\Cal Q(j))\ne\emptyset$. 
Then there is an immersion $h\in \Cal W(i)\cap \Cal Q(j)$, so there 
exists a unique $f_0\in\Diff(M)$ (given by $f_0=p\o\ph_i(h)$) such 
that $h\o f_0\i\in\Cal Q(i)$. If we consider $j\o f_0\i$ instead of 
$j$ and call it again $j$, we have 
$\Cal Q(i)\cap\Cal Q(j)\ne\emptyset$ and consequently 
$\Cal U(i)\cap\Cal U(j)\ne\emptyset$.
Then the chart change is given as follows:
$$\gather
\ph_i\o(\pi|\Cal Q(i))\i\o\pi\o(\tau^j)_*:
     C^\infty_c(U^j)\to C^\infty_c(U^i)\\
s\mapsto \tau^j\o s 
     \mapsto \ph_i(\tau^j\o s)\o(p^i\o\ph_i(\tau^j\o s))\i.
\endgather$$
This is of the form $s\mapsto \be\o s$ for a locally defined 
diffeomorphism $\be:\Cal N(j)\to\Cal N(i)$ which is not fiber 
respecting, followed by $h\mapsto h\o(p^i\o h)\i$. Both composants 
are smooth by the general properties of manifolds of mappings. 
So the chart change is smooth.

We have to show that the quotient space 
$\Imm_{\text{prop,free}}(M,N)/\Diff(M)$ is Hausdorff. This will be done in 
section \nmb!{2} below.
\qed\enddemo

\heading\totoc\nmb0{2}. Some orbit spaces are Hausdorff\endheading

\proclaim{\nmb.{2.1}. Theorem} The orbit space
$\Imm_{\text{prop}}(M,N)/\Diff(M)$ of the space of all proper
immersions under the action of the diffeomorphism group is
Hausdorff in the quotient topology.
\endproclaim

The proof will occupy the rest of this section. We want to point out
that we believe that the whole orbit space 
$\Imm(M,N)/\Diff(M)$ is Hausdorff, but that we were unable to prove
this.

\proclaim{\nmb.{2.2}. Lemma} Let $i$ and $j\in \Imm_{\text{prop}}(M,N)$ with 
$i(M)\ne j(M)$ in $N$. Then their projections $\pi(i)$ and $\pi(j)$ 
are different and can be separated by open subsets in 
$\Imm_{\text{prop}}(M,N)/\Diff(M)$.
\endproclaim

\demo{Proof} We suppose that $i(M)\nsubseteq \overline{j(M)} = j(M)$ 
(since proper immersions have closed images). 
Let $y_0\in i(M)\setminus\overline{j(M)}$, then we choose open 
neighborhoods $V$ of $y_0$ in $N$ and $W$ of $j(M)$ in $N$ such that 
$V\cap W = \emptyset$. We consider the sets
$$\align
\Cal V &:= \{k\in \Imm_{\text{prop}}(M,N):k(M)\cap V\neq\emptyset\}\quad\text{and}\\
\Cal W &:= \{k\in \Imm_{\text{prop}}(M,N):k(M)\subseteq W\}.
\endalign$$
Then $\Cal V$ and $\Cal W$ are $\Diff(M)$-saturated 
disjoint open neighborhoods of $i$ and $j$, respectively, so 
$\pi(\Cal V)$ and $\pi(\Cal W)$ separate $\pi(i)$ and $\pi(j)$ in 
$\Imm_{\text{prop}}(M,N)/\Diff(M)$.
\qed\enddemo

\subheading{\nmb.{2.3}} For a proper immersion $i:M\to N$ and 
$x\in i(M)$ let $\de(x)\in \Bbb N$ be the number of points in 
$i\i(x)$. Then $\de:i(M)\to \Bbb N$ is a mapping.

\proclaim{Lemma} The mapping $\de:i(M)\to \Bbb N$ is upper 
semicontinuous, i.e. $\{x\in i(M):\de(x)\leq k\}$ is open in 
$i(M)$ for each $k$.
\endproclaim

\demo{Proof}
Let $x\in i(M)$ with $\de(x)=k$ and let $i\i(x)=\{\row y1k\})$. Then 
there are pairwise disjoint open neighborhoods $W_n$ of $y_n$ in $M$ 
such that $i|W_n$ is an embedding for each $n$. The set 
$M\setminus(\bigcup_n W_n)$ is closed in $M$, and since $i$ is proper 
the set $i(M\setminus(\bigcup_n W_n))$ is also closed in $i(M)$ and 
does not contain $x$. So there is an open neighborhood $U$ of $x$ in 
$i(M)$ which does not meet $i(M\setminus(\bigcup_n W_n))$. Then 
obviously $\de(z)\leq k$ for all $z\in U$.
\qed\enddemo

\subheading{\nmb.{2.4}} We consider two proper immersions $i_1$ 
and $i_2\in \Imm_{\text{prop}}(M,N)$ such that 
$i_1(M)=i_2(M)=:L\subseteq N$. 
Then we have mappings $\de_1, \de_2: L\to \Bbb N$ as in \nmb!{2.3}.

\proclaim{\nmb.{2.5}. Lemma} In the situation of \nmb!{2.4}, if 
$\de_1\ne\de_2$ then the projections $\pi(i_1)$ and $\pi(i_2)$ are 
different and can be 
separated by disjoint open neighborhoods in 
$\Imm_{\text{prop}}(M,N)/\Diff(M)$.
\endproclaim

\demo{Proof}
Let us suppose that $m_1=\de_1(y_0)\ne\de_2(y_0)=m_2$. 
There is a small connected open neighborhood $V$ of $y_0$ in $N$ such 
that $i_1\i(V)$ has $m_1$ connected components and $i_2\i(V)$ has 
$m_2$ connected components. This assertions 
describe Whitney $C^0$-open neighborhoods in $\Imm_{\text{prop}}(M,N)$ of 
$i_1$ and $i_2$ which are closed under the action of $\Diff(M)$, 
respectively. Obviously these two neighborhoods are disjoint.
\qed\enddemo

\subheading{\nmb.{2.6} } We assume now for the rest of this section 
that we are given two immersions $i_1$ and $i_2\in \Imm_{\text{prop}}(M,N)$ with
$i_1(M)=i_2(M)=:L$ such that the functions from \nmb!{2.4} are equal:
$\de_1=\de_2=:\de$.

Let $(L_\be)_{\be\in B}$ be the partition of $L$ consisting of all 
pathwise connected components of level sets $\{x\in L:\de(x)=c\}$, $c$ 
some constant. 

Let $B_0$ denote the set of all $\be\in B$ such that the interior of 
$L_\be$ in $L$ is not empty. Since $M$ is second countable, $B_0$ is
countable. \newline
{\bf Claim.} $\bigcup_{\be\in B_0}L_\be$ is dense in $L$. \newline
Let $k_1$ be the smallest number in $\de(L)$ and let $B_1$ be the set 
of all $\be\in B$ such that $\de(L_\be)=k_1$. Then by lemma 
\nmb!{2.3} each $L_\be$ for $\be\in B_1$ is open. Let $L^1$ be the 
closure of $\bigcup_{\be\in B_1}L_\be$. Let $k_2$ be the smallest 
number in $\de(L\setminus L^1)$ and let $B_2$ be the set of all 
$\be\in B$ with $\be (L_\be)=k_2$ and 
$L_\be\cap(L\setminus L^1)\ne\emptyset$. Then by lemma \nmb!{2.3} again  
$L_\be\cap(L\setminus L^1)\ne\emptyset$ is open in $L$ so $L_\be$ has 
non empty interior for each $\be\in B_2$. Then let $L^2$ denote the 
closure of $\bigcup_{\be\in B_1\cup B_2}L_\be$ and continue the process. 
Since by lemma \nmb!{2.3} we always find new $L_\be$ with non empty 
interior, we finally exhaust $L$ and the claim follows.

Let $(M^1_\la)_{\la\in C^1}$ be a suitably chosen cover of $M$ by subsets 
of the sets $i_1\i(L_\be)$ such that each 
$i_2| \oin M^1_\la$ is an embedding for each 
$\la$. Let $C^1_0$ be the set of all $\la$ such that $M^1_\la$ has non 
empty interior. Let similarly $(M^2_\mu)_{\mu\in C^2}$ be a cover for $i_2$.
Then there are at most countably many sets $M^1_\la$ with  
$\la\in C^1_0$, the union
$\bigcup_{\la\in C^1_0}\oin M^1_\la$ is dense and 
consequently $\bigcup_{\la\in C^1_0}\overline{M^1_\la}= M$; similarly 
for the $M^2_\mu$.

\subheading{\nmb.{2.7}. Procedure} Given immersions $i_1$ and $i_2$ 
as in \nmb!{2.6} we will try to construct a diffeomorphism $f:M\to M$ 
with $i_2\o f=i_1$. If we meet an obstacle to the construction this 
will give us enough control on the situation to separate $i_1$ and $i_2$.

Choose $\la_0\in C^1_0$ so that 
$\oin M^1_{\la_0}\ne\emptyset$. Then 
$i_1:\oin M^1_{\la_0}\to L_{\be_1(\la_0)}$ is an 
embedding, where $\be_1:C^1\to B$ is the mapping satisfying 
$i_1(M^1_\la)\subseteq L_{\be_1(\la)}$ for all $\la\in C^1$.

Now we choose $\mu_0\in \be_2\i\be_1(\la_0)\subset C^2_0$ such that 
$f:=(i_2|\oin M^2_{\mu_0})\i\o 
     i_1|\oin M^1_{\la_0}$
is a diffeomorphism $\oin M^1_{\la_0}\to \oin M^2_{\mu_0}$.
Note that $f$ is uniquely determined by the choice of $\mu_0$, if it 
exists, by lemma \nmb!{1.3}. So we will repeat the following 
construction for every $\mu_0\in \be_2\i\be_1(\la_0)\subset C^2_0$.

Now we try to extend $f$. 
We choose $\la_1\in C^1_0$ such that  
$\overline M^1_{\la_0}\cap\overline M^1_{\la_1}\ne \emptyset$.

{\bf Case a.} Only $\la_1=\la_0$ is possible, so $M^1_{\la_0}$ 
is dense in $M$ 
since $M$ is connected and we may extend $f$ by continuity to a 
diffeomorphism $f:M\to M$ with $i_2\o f=i_1$.

{\bf Case b.} We can find $\la_1\ne\la_0$. We choose 
$x\in\overline M^1_{\la_0}\cap\overline M^1_{\la_1}$ and a sequence 
$(x_n)$ in $M^1_{\la_0}$ with $x_n\to x$. Then we have a sequence
$(f(x_n))$ in $B$.

{\bf Case ba.} $y:=\lim f(x_n)$ exists in $M$. Then there is 
$\mu_1\in C^2_0$ such that 
$y\in \overline M^2_{\mu_0}\cap\overline M^2_{\mu_1}$.

Let $U^1_{\al_1}$ be an open neighborhood of $x$ in $M$ such that 
$i_1| U^1_{\al_1}$ is an embedding and let similarly 
$U^2_{\al_2}$ be an open neighborhood of $y$ in $M$ such that
$i_2| U^2_{\al_2}$ is an embedding. We consider now the 
set $i_2\i i_1(U^1_{\al_1})$. There are two cases possible.

{\bf Case baa.} The set $i_2\i i_1(U^1_{\al_1})$ is a neighborhood of 
$y$. Then we extend $f$ to $i_1\i(i_1(U^1_{\al_1})\cap i_2(U^2_{\al_2}))$ 
by $i_2\i\o i_1$. Then $f$ is defined on some open subset of 
$\oin M^1_{\la_1}$ and by the situation chosen in \nmb!{2.6} $f$ extends 
to the whole of $\oin M^1_{\la_1}$.

{\bf Case bab.} The set $i_2\i i_1(U^1_{\al_1})$ is not a neighborhood of 
$y$. This is a definite obstruction to the extension of $f$.

{\bf Case bb.} The sequence $(x_n)$ has no limit in $M$. This is a 
definite obstruction to the extension of $f$.

If we meet an obstruction we stop and try another $\mu_0$. If for all 
admissible $\mu_0$ we meet obstructions we stop and remember the data. 
If we do not 
meet an obstruction we repeat the construction with some obvious 
changes.

\proclaim{\nmb.{2.8}. Lemma} The construction of \nmb!{2.7} in the 
setting of \nmb!{2.6} either produces a diffeomorphism $f:M\to M$ 
with $i_2\o f=i_1$ or we may separate $i_1$ and $i_2$ by open sets in 
$\Imm_{\text{prop}}(M,N)$ which are saturated with respect to the action of 
$\Diff(M)$
\endproclaim

\demo{Proof}
If for some $\mu_0$ we do not meet any obstruction in the 
construction \nmb!{2.7}, the resulting $f$ is defined on the whole of 
$M$ and it is a 
continuous mapping $M\to M$ with $i_2\o f=i_1$. Since $i_1$ and $i_2$ 
are locally embeddings, $f$ is smooth and of maximal rank. Since $i_1$ 
and $i_2$ are proper, $f$ is proper. So the image of $f$ is open and 
closed and since $M$ is connected, $f$ is a surjective local 
diffeomorphism, thus a covering mapping $M\to M$. But since 
$\de_1=\de_2$ the mapping $f$ must be a 1-fold covering, so a 
diffeomorphism.

If for all $\mu_0\in \be_2\i\be_1(\la_0)\subset C^2_0$ we meet 
obstructions we choose small mutually distinct open neighborhoods 
$V^1_\la$ of the sets $i_1(M^1_\la)$. We consider the Whitney 
$C^0$-open neighborhood $\Cal V_1$ of $i_1$  consisting of all 
immersions $j_1$ with $j_1(M^1_\la)\subset V^1_\la$ for all $\la$.
Let $\Cal V_2$ be a similar neighborhood of $i_2$.

We claim that $\Cal V_1 \o \Diff(M)$ and $\Cal V_2\o\Diff(M)$ are 
disjoint. For that it suffices to show that for any $j_1\in\Cal V_1$ 
and $j_2\in \Cal V_2$ there does not exist a diffeomorphism 
$f\in\Diff(M)$ with $j_2\o f=j_1$. For that to be possible the 
immersions $j_1$ and $j_2$ must have the same image $L$ and the same 
functions $\de(j_1)$, $\de(j_2):L\to \Bbb N$. But now the 
combinatorial relations of the slightly distinct new sets $M^1_\la$, 
$L_\be$, and $M^2_\mu$ are contained in the old ones, so any try to 
construct such a diffeomorphism $f$ starting from the same $\la_0$ 
meets the same obstructions.
\qed\enddemo

\heading\totoc\nmb0{3}. Singular orbits \endheading

\subheading{\nmb.{3.1}} Let $i\in\Imm(M,N)$ be an immersion which is 
not free. Then we have a nontrivial isotropy subgroup 
$\Diff_i(M)\subset \Diff(M)$ consisting of all $f\in\Diff(M)$ with 
$i\o f= i$.

\proclaim{Lemma} Then the isotropy subgroup $\Diff_i(M)$ acts 
properly discontinuously on $M$, so the projection 
$q_1:M\to M_1:= M/\Diff_i(M)$ is a covering map and a submersion for 
a unique structure of a smooth manifold on $M_1$. There is an 
immersion $i_1:M_1\to N$ with $i=i_1\o q_1$. In particular $\Diff_i(M)$ 
is countable, and finite if $M$ is compact.
\endproclaim

\demo{Proof}
We have to show that for each $x\in M$ there is an open neighborhood 
$U$ such that $f(U)\cap U=\emptyset$ for 
$f\in \Diff_i(M)\setminus \{Id\}$. We consider the setup \nmb!{1.1} for 
$i$. By the proof of \nmb!{1.3} we have 
$f(U^i_\al)\cap U^i_\al=\{x\in U^i_\al: f(x)=x\}$ for any 
$f\in\Diff_i(M)$. If $f$ has a fixed point then by \nmb!{1.3} $f=Id$, 
so $f(U^i_\al)\cap U^i_\al=\emptyset$ for all 
$f\in\Diff_i(M)\setminus\{Id\}$. The rest is clear.
\qed\enddemo
The factorized immersion $i_1$ is in general not a free immersion. 
The following is an example for that: Let 
$$M_0 @>>\al> M_1 @>>\be> M_2 @>>\ga> M_3$$
be a sequence of covering maps with fundamental groups 
$1\to G_1 \to G_2 \to G_3$. Then the group of deck transformations of 
$\ga$ is given by $\Cal N_{G_3}(G_2)/G_2$, the normalizer of $G_2$ in 
$G_3$, and the group of deck transformations of 
$\ga\o\be$ is $\Cal N_{G_3}(G_1)/G_1$. We can 
easily arrange that $\Cal N_{G_3}(G_2)\nsubseteq\Cal N_{G_3}(G_1)$, 
then $\ga$ admits deck transformations which do not lift to $M_1$. 
Then we thicken all spaces to manifolds, so that $\ga\o\be$ plays the 
role of the immersion $i$. 

\proclaim{\nmb.{3.2}. Theorem} Let $i\in\Imm(M,N)$ be an immersion 
which is not free. Then there is a covering map $q_2:M\to M_2$ 
which is also a submersion such that $i$ factors to an immersion 
$i_2:M_2\to N$ which is free.
\endproclaim

\demo{Proof}
Let $q_0:M_0\to M$ be the universal covering of $M$ and 
consider the immersion $i_0=i\o q_0:M_0\to N$ and 
its isotropy group $\Diff_{i_0}(M_0)$. By \nmb!{3.1} it 
acts properly discontinuously on $M_0$ and we have a submersive 
covering $q_{02}:M_0\to M_2$ and an immersion 
$i_2:M_2\to N$ with 
$i_2\o q_{02}=i_0=i\o q_0$. By comparing the respective groups of 
deck transformations it is easily seen that 
$q_{02}:M_0\to M_2$ factors over
$q_1\o q_0:M_0\to M\to M_1$ 
to a covering $q_{12}:M_1\to M_2$. The mapping 
$q_2:=q_{12}\o q_1:M\to M_2$ is the looked for covering: If 
$f\in \Diff(M_2)$ fixes $i_2$, it lifts to a 
diffeomorphism $f_0\in \Diff(M_0)$ which fixes $i_0$, 
so is in $\Diff_{i_0}(M_0)$, so $f=Id$.
\qed\enddemo

\subheading{\nmb.{3.3}. Convention} In order to avoid complications 
we assume that from now on $M$ is such a manifold that 
\roster
\item For any covering $M\to M_1$, any diffeomorphism $M_1\to M_1$ 
     admits a lift $M\to M$.
\endroster
If $M$ is simply connected, condition \therosteritem1 is satisfied. 
Also for $M=S^1$ condition \therosteritem1 is easily seen to be 
valid. So what follows is applicable to loop spaces.

Condition \therosteritem1 implies that in the proof of 
\nmb!{3.2} we have $M_1=M_2$.

\subheading{\nmb.{3.4}. Description of a neighborhood of a singular 
orbit} 
Let $M$ be a manifold satisfying \nmb!{3.3}.\therosteritem1.
In the situation of \nmb!{3.1} we consider the normal bundles 
$p_i:\Cal N(i)\to M$ and 
$p_{i_1}:\Cal N(i_1)\to M_1$. Then the 
covering map $q_1:M\to M_1$ lifts uniquely to a vector bundle 
homomorphism $\Cal N(q_1):\Cal N(i)\to \Cal N(i_1)$ which is 
also a covering map, such that $\tau^{i_1}\o\Cal N(q_1)=\tau^i$.

We have $M_1= M/\Diff_i(M)$ and the group $\Diff_i(M)$ acts also as 
the group of deck transformations  of the covering 
$\Cal N(q_1):\Cal N(i)\to \Cal N(i_1)$ by 
$\Diff_i(M)\ni f\mapsto \Cal N(f)$, where 
$$\CD
\Cal N(i) @>>\Cal N(f)>  \Cal N(i)\\
@VVV                      @VVV\\
M          @>>f>          M
\endCD$$
is a vector bundle isomorphism for each $f\in \Diff_i(M)$.
If we equip $\Cal N(i)$ and $\Cal N(i_1)$ with the fiber Riemann 
metrics induced from the fixed Riemannian metric $g$ on $N$, the 
mappings $\Cal N(q_1)$ and all $\Cal N(f)$ are fiberwise linear 
isometries.

Let us now consider the right action of $\Diff_i(M)$ on the space of 
sections $C^\infty_c(\Cal N(i))$ given by 
$f^*s:=\Cal N(f)\i\o s\o f$. 

 From the proof of theorem \nmb!{1.5} we recall now the sets
$$\CD 
C^\infty(M,\Cal N(i))\supset \Cal V(i) @<<\ph_i< \Cal U(i)\\
@AAA                                              @AAA\\
C^\infty_c(\Cal N(i))\supset C^\infty_c(U^i) @<\ph_i<< \Cal Q(i).
\endCD$$
All horizontal mappings are again diffeomorphisms and the vertical 
mappings are inclusions. But since the action of $\Diff(M)$ on $i$ is 
not free we cannot extend the splitting submanifold $\Cal Q(i)$ to an 
orbit cylinder as we did in the proof on theorem \nmb!{1.5}. 
$\Cal Q(i)$ is again a smooth transversal for the orbit though $i$.

For any $f\in \Diff(M)$ and 
$s\in C^\infty_c(U^i)\subset C^\infty_c(\Cal N(i))$ we have 
$$\ph_i\i(f^*s)= \tau^i_*(f^*s) = \tau^i_*(s)\o f.$$
So the space $q_1^*C^\infty_c(\Cal N(i_1))$ of all sections of 
$\Cal N(i)\to M$
which factor to sections of $\Cal N(i_1)\to M_1$, is exactly the 
space of all fixed points of the action of $\Diff_i(M)$ on 
$C^\infty_c(\Cal N(i))$; and they are mapped by $\tau^i_*=\ph_i\i$ to 
immersions in $\Cal Q(i)$ which have again $\Diff_i(M)$ as isotropy 
group.

If $s\in C^\infty_c(U^i)\subset C^\infty_c(\Cal N(i))$ is an 
arbitrary  section, the orbit through $\tau^i_*(s)\in \Cal Q(i)$ hits 
the transversal $\Cal Q(i)$ again in the points $\tau^i_*(f^*s)$ for 
$f\in\Diff_i(M)$.

We summarize all this in the following theorem:

\proclaim{\nmb.{3.5}. Theorem} Let $M$ be a manifold satisfying 
condition \therosteritem1 of \nmb!{3.3}. Let $i\in Imm(M,N)$ be an 
immersion which is not free, i.e. has non trivial isotropy group 
$\Diff_i(M)$. 

Then in the setting and notation of \nmb!{3.4}
in the following commutative diagram the bottom mapping
$$\CD
\Imm_{\text{free}}(M_1,N) @>(q_1)^*>> \Imm(M,N)\\
@V\pi VV                               @VV\pi V \\
\Imm_{\text{free}}(M_1,N)/\Diff(M_1) @>>>\Imm(M,N)/\Diff(M)
\endCD$$ 
is the inclusion of a (possibly non Hausdorff) manifold, 
the stratum of $\pi(i)$ in the 
stratification of the orbit space. This stratum consists of the 
orbits of all immersions which have $\Diff_i(M)$ as isotropy group.
\endproclaim

\subheading{\nmb.{3.6}. The orbit structure} 
We have the following description of the orbit structure near $i$ 
in $\Imm(M,N)$: For fixed $f\in \Diff_i(M)$ the set of fixed points
$\operatorname{Fix}(f):=\{j\in\Cal Q(i):j\o f=j\}$ is called a 
\idx{\it generalized wall}. The union of all generalized walls is called the 
\idx{\it diagram} $\Cal D(i)$ of $i$. A connected component of the complement
$\Cal Q(i)\setminus \Cal D(i)$ is called a \idx{\it generalized Weyl
chamber}. 
The group $\Diff_i(M)$ maps walls to walls and chambers to chambers.
The immersion $i$ lies in every wall. 

We shall see shortly that there is only one chamber and that the 
situation is rather distinct from that of reflection groups.

If we view the diagram in the space
$C^\infty_c(U^i)\subset C^\infty_c(\Cal N(i))$ which is diffeomorphic 
to $\Cal Q(i)$, then it consist of
traces of closed linear subspaces, because
the action of $\Diff_i(M)$ on $C^\infty_c(\Cal N(i))$ consists of
linear isometries in the following way. Let us tensor the vector 
bundle $\Cal N(i)\to M$ with the natural line bundle of half 
densities on $M$, and let us remember one positive half density to
fix an isomorphism with the original bundle. 
Then $\Diff_i(M)$ still acts on this new bundle $\Cal N_{1/2}(i)\to M$ 
and the pullback
action on sections with compact support is isometric for the inner product
$$\langle s_1,s_2 \rangle := \int_M g(s_1,s_2).$$
We consider the walls and chambers now extended to the whole space in 
the obvious manner. 

\proclaim{\nmb.{3.7}. Lemma} Each wall in $C^\infty_c(\Cal N_{1/2}(i))$ 
is a closed linear subspace of infinite codimension. Since there are 
at most countably many walls, there is only one chamber. 
\endproclaim
\demo{Proof}
 From the proof of lemma \nmb!{3.1} we know that 
$f(U^i_\al)\cap U^i_\al=\emptyset$ for all $f\in\Diff_i(M)$ and all 
sets $U^i_\al$ from the setup \nmb!{1.1}. Take a section $s$ in the 
wall of fixed points of $f$.
Choose a section $s_\al$ with support in some $U^i_\al$ and let 
the section $s$ be defined by 
$s|U^i_\al=s_\al|U^i_\al$, $s|f\i(U^i_\al)=-f^*s_\al$, 0 elsewhere.
Then obviously $\langle s,s'\rangle=0$ for all $s'$ in the wall of 
$f$. But this construction furnishes an infinite dimensional space 
contained in the orthogonal complement of the wall of $f$.
\qed\enddemo


\Refs

\ref 
\by Binz, Ernst; Fischer, Hans R.
\paper The manifold of embeddings of a closed manifold
\inbook Proc. Differential geometric methods in theoretical physics, Clausthal 1978
\publ Springer Lecture Notes in Physics 139 
\yr 1981 
\endref 

\ref \by Fr\"olicher, A.; Kriegl, A. \book Linear
spaces and differentiation theory \bookinfo Pure and Applied
Mathematics \publ J. Wiley \publaddr Chichester \yr 1988 \endref

\ref  
\by Kriegl, A.; Michor, P. W.  
\paper A convenient setting for real analytic mappings 
\paperinfo 52 p., to appear  
\jour Acta Mathematica 
\yr 1990 
\endref

\ref 
\by Michor, P. W.
\paper  Manifolds of smooth maps 
\jour Cahiers Topol. Geo. Diff.
\vol 19 
\yr 1978
\pages 47--78
\endref

\ref 
\by Michor, P. W.
\paper Manifolds of smooth maps II: The Lie group of diffeomorphisms of a non compact smooth manifold
\jour Cahiers Topol. Geo. Diff. 
\vol 21 
\yr 1980a
\pages 63--86
\endref

\ref 
\by Michor, P. W.
\paper Manifolds of smooth maps III: The principal bundle of embeddings of a non compact smooth manifold
\jour Cahiers Topol. Geo. Diff. 
\vol 21 
\yr 1980b
\pages 325--337
\endref

\ref\by Michor, P. W. \book Manifolds of differentiable mappings
\publ Shiva \yr 1980c \publaddr Orpington \endref

\ref \by Michor, P. W. \paper Manifolds of smooth
mappings IV: Theorem of De~Rham \jour Cahiers Top. Geo. Diff. 
\vol 24 \yr 1983 \pages 57--86 \endref

\ref \by Michor, P. W. \paper Gauge theory for
diffeomorphism groups \inbook Proceedings of the Conference on
Differential Geometric Methods in Theoretical Physics, Como
1987, K. Bleuler and M. Werner (eds.)
\publ Kluwer \publaddr Dordrecht \yr 1988 \pages 345--371 \endref 

\ref 
\by Weinstein, Alan
\paper Symplectic manifolds and their Lagrangian manifolds
\jour Advances in Math.
\vol 6
\pages 329--345
\yr 1971
\endref

\endRefs
\enddocument